\RequirePackage{ifpdf}
\ifpdf 
\documentclass[pdftex]{sigma}
\else
\documentclass{sigma}
\fi

\begin{document}

\allowdisplaybreaks

\renewcommand{\PaperNumber}{017}

\FirstPageHeading

\ShortArticleName{Orthogonality Relations for Multivariate Krawtchouk Polynomials}

\ArticleName{Orthogonality Relations\\ for Multivariate Krawtchouk Polynomials}

\Author{Hiroshi MIZUKAWA}

\AuthorNameForHeading{H.~Mizukawa}

\Address{Department of Mathematics, National Defense Academy of Japan,
Yokosuka 239-8686, Japan}

\Email{\href{mailto:mzh@nda.ac.jp}{mzh@nda.ac.jp}}

\ArticleDates{Received September 08, 2010, in f\/inal form February 18, 2011;  Published online February 22, 2011}

\Abstract{The orthogonality relations of multivariate Krawtchouk polynomials are discussed.
In case of  two variables, the necessary and suf\/f\/icient conditions of orthogonality
is given by Gr\"{u}nbaum  and Rahman in~[{\it SIGMA} {\bf 6} (2010), 090, 12~pages].
In this study,  a~simple proof of the necessary and suf\/f\/icient condition of orthogonality
is given for a general case.}

\Keywords{multivariate orthogonal polynomial; hypergeometric function}
\Classification{33C45}

\section{Introduction}
Consider
\[
X(n,N)=\{\boldsymbol{x}= (x_{0},x_{1},\dots,x_{n-1})\in {\mathbb N_{0}^{n}} \mid |\boldsymbol{x}|=N\},
\]
where ${\mathbb N_0}$ is the set of nonnegative integers and $|\boldsymbol{x}|=x_{0}+x_{1}+\cdots+x_{n-1}$.
We recall the multinomial coef\/f\/icient
\[
\binom{N}{\boldsymbol{x}}=\binom{N}{x_{0},\dots,x_{n-1}}
=(-1)^{x_{1}+\cdots+x_{n-1}}\frac{(-N)_{x_{1}+\cdots+x_{n-1}}}{x_{1}!\cdots x_{n-1}!}
\] for $\boldsymbol{x} \in X(n,N)$.
Let $M_{n}(R)$ be the set of all $n \times n$ matrices over a set $R$.
We f\/ix $\boldsymbol{x} \in X(n,N)$ and $A=(a_{ij})_{1\leq i,j\leq n-1} \in M_{n-1}({\mathbb C})$.
We def\/ine the functions $\phi_A(\boldsymbol{x};\boldsymbol{m})$ of $\boldsymbol{m}=(m_{0},\dots,m_{n-1}) \in X(n,N)$ by the following generating function
\begin{gather}\label{gen}
\Phi_{N}(A;\boldsymbol{x})=
\prod_{i=0}^{n-1}\left(\sum_{j=0}^{n-1}a_{ij}t_{j}\right)^{x_{i}}=\sum_{\boldsymbol{m} \in X(n,N)}
\binom{N}{\boldsymbol{m}}\phi_A(\boldsymbol{x};\boldsymbol{m})\boldsymbol{t}^{\boldsymbol{m}},
\end{gather}
where $\boldsymbol{t}^{\boldsymbol{m}}=t_{0}^{m_{0}}t_{1}^{m_{1}}\cdots t_{n-1}^{m_{n-1}}$
and $a_{0j}=a_{i0}=1$ for $0\leq i,j\leq n-1$.
We know a hypergeometric expression of $\phi_A(\boldsymbol{x};\boldsymbol{m})$
\begin{gather}\label{hyp}
\phi_A(\boldsymbol{x};\boldsymbol{m})
=\sum_{\substack{\sum_{i,j}c_{ij}\leq N\\
(c_{ij})\in M_{n-1}({\mathbb N_{0}})}}
\frac{\prod\limits_{i=1}^{n-1}(-x_{i})_{\sum\limits_{j=1}^{n-1}c_{ij}}
\prod\limits_{j=1}^{n-1}(-m_{j})_{\sum\limits_{i=1}^{n-1}c_{ij}}}
{(-N)_{\sum_{i,j}c_{ij}}} \; \frac{\prod
(1-a_{ij})^{c_{ij}}}{\prod c_{ij}!},
\end{gather}
where $M_{n-1}({\mathbb N}_{0})$ is the set of square matrices of degree $n-1$ with nonnegative integer
elements.
We prove the formula (\ref{hyp}) in the last section of this paper.

This type of hypergeometric functions was originally def\/ined by Aomoto and Gel'fand
for general parameters.
We are interested in the aspects of discrete orthogonal polynomials of these functions
with weights
\[
b_{n}(\boldsymbol{x};N;\boldsymbol{\eta}_{(i)})=\binom{N}{\boldsymbol{x}}\prod_{j=0}^{n-1}\eta_{ji}^{x_{j}}\]
for $\boldsymbol{\eta}_{(i)}=(\eta_{0i},\dots,\eta_{n-1i}) \in {\mathbb C^*}^n$  $(i=1,2)$,
where ${\mathbb C}^*={\mathbb C}\backslash \{0\}$.
For a special case, when $n=2$, they are well known and called the Krawtchouk polynomials.
For general values of $n$, the author shows that such
 orthogonal polynomials  appear as the zonal spherical functions of Gel'fand pairs of complex
ref\/lection groups~\cite{m}.
In general,  the author and H.~Tanaka give the orthogonality relation of $\phi_A(\boldsymbol{x};\boldsymbol{m})$s
by using the character algebras~\cite{mt}.
R.C.~Grif\/f\/iths shows that the polynomials def\/ined by the generating function (\ref{gen}) are mutually orthogonal \cite{gr}.
In their paper~\cite{rg}, Gr\"{u}nbaum and Rahman discuss and determine the necessary and suf\/f\/icient
conditions of the orthogonality of $\phi_A(\boldsymbol{x};\boldsymbol{m})$s for $A \in M_{2}({\mathbb C})$,
which are proved by analytic methods.
In these four literature, the authors consider the case that weights $b_{n}(\boldsymbol{x};N;\boldsymbol{\eta}_{(i)})$
are positive.

In this study, we give a linear algebraic proof of
the Gr\"{u}nbaum and Rahman's condition for general values of $n$ and
for arbitrary weights including the complex case.
Our proof of the suf\/f\/icient condition is close to~\cite{gr} and~\cite{mt}.

For $A \in M_{n-1}({\mathbb C})$, we def\/ine a matrix $A_{0}=(a_{ij})_{0\leq i,j\leq n-1}
 \in M_{n}({\mathbb C})$ by
substituting $a_{0j}=a_{i0}=1$ ($0 \leq \forall\,  i,j \leq n-1$).
Our main result is as follows:

\begin{theorem}\label{a1}
The following are equivalent.
\begin{enumerate}\itemsep=0pt
\item[$(a)$]
The orthogonality relation
\[\sum_{\boldsymbol{x} \in X(n,N)}b_{n}(\boldsymbol{x};N;\boldsymbol{\eta}_{(1)})
\phi_A(\boldsymbol{x};\boldsymbol{m})\overline{\phi_A(\boldsymbol{x};\boldsymbol{m}')}=\delta_{\boldsymbol{m},\boldsymbol{m}'}
\frac{\boldsymbol{\eta}_{(2)}^{\boldsymbol{m}}}{\binom{N}{\boldsymbol{m}}}
\]
holds for $\forall \, \boldsymbol{m},\boldsymbol{m}' \in X(n,N)$.
\item[$(b)$]
A relation
\begin{gather}\label{ada}
A_{0}^{*}D_{1}{A_{0}}= \zeta D_{2}
\end{gather}
holds for some $N$th root of unity $\zeta$. Here $A^*_{0}$ is the conjugate transpose of $A_{0}$ and
$D_{i}={\rm diag}(\eta_{0i},\eta_{1i},\dots,\eta_{n-1i})\in {\rm GL}_{n}({\mathbb C})$  is a  diagonal matrix $(i=1,2)$.
\end{enumerate}
\end{theorem}
 \begin{remark}
We assume that  the diagonal elements
of $D_{1}$ and $D_{2}$ appearing in the above-mentioned theorem are real.
Thus, one can recover the  formula (1.18) of Gr\"{u}nbaum and Rahman's paper \cite{rg}
by substituting $n=3$ and
 $A=
 \left[\begin{array}{cc}
 1-u_{1}&1-u_{2}\\
 1-v_{1}&1-v_{2}
 \end{array}
 \right] \in M_{2}({\mathbb C})$.
In this case, since  the diagonal elements of $A_{0}^{*}D_{1}{A_{0}}$ and $D_{2}$ are positive, $\zeta=1$.
\end{remark}

\begin{remark}[\cite{mt}]
We assume that a pair of f\/inite groups $(G,H)$ is a Gel'fand pair and $A_{0}$
is the table of the zonal spherical functions of $(G,H)$.
Let $D_{0},\dots,D_{n-1}$ be the double cosets of $H$ in $G$, and
$d_{0},\dots,d_{n}$ be the dimensions of the  irreducible components  of $1_{H}^G$.
Put $D_{1}={\rm diag}(|D_{0}|,\dots,|D_{n-1}|)$
and $D_{2}={\rm diag}(|G|/d_{0},\dots,|G|/d_{n-1})$.
Then~(\ref{ada}) holds from the orthogonality relation of the zonal spherical functions.
Furthermore $\phi_A(\boldsymbol{x};\boldsymbol{m})$'s are realized as the zonal spherical functions of a Gel'fand pair
 $(G \wr S_{N},H \wr S_{N})$. Therefore they satisfy the orthogonality relation (a) in the theorem.
In general, $A_{0}$ is  an eigenmatrix of a character algebra is considered in~\cite{mt}.
 \end{remark}

 This paper organized as follows.
First, we prove the main theorem in the next section.
Second, we prove (\ref{hyp}) in the last section. It seems to be the f\/irst explicit proof of this fact.

\section{Proof of Theorem \ref{a1}}\label{section2}

Let ${\bf e}_{i}=(\delta_{0i},\delta_{1i},\dots,\delta_{ni})$ be an $i$th unit vector ($0 \leq i \leq n-1$).
We put ${\bf s}=(s_{0},s_{1},\dots,s_{n-1})$ and ${\bf t}
=(t_{0},t_{1},\dots,t_{n-1})$. We compute
\begin{gather}\label{21}
\prod_{i=0}^{n-1}\big(s_{i}{\bf e}_{i} A_{0}\, {}^t {\bf t}\big)^{x_{i}}=
\prod_{i=0}^{n-1}\left(\sum_{j=0}^{n-1}a_{ij}t_{j}\right)^{x_{i}}
\prod_{i=0}^{n-1}s_{i}^{x_{i}}=\Phi_{N}(A,\boldsymbol{x})\prod_{i=0}^{n-1}s_{i}^{x_{i}}.
\end{gather}
By multiplying  (\ref{21}) with multinomial coef\/f\/icient, we have
\begin{gather}\label{pdt}
\sum_{\boldsymbol{x} \in X(n,N)}\binom{N}{\boldsymbol{x}}\prod_{i=0}^{n-1}\big(s_{i}{\bf e}_{i} A_{0} {\, {}^t\bf t}\big)^{x_{i}}
=\left(\sum_{i=0}^{n-1}s_{i}{\bf e}_{i} A_{0} {\, {}^t\bf t}\right)^N
=\big({\bf s}A_{0}{\, {}^t\bf t}\big)^N.
\end{gather}

First, we assume  (\ref{ada}), and then
change the variables, say
\[{\bf s}={\bf u}A_{0}^{*}D_{1},\]
where  ${\bf u}=(u_{0},u_{1},\dots,u_{n-1})$.
This change of variables is same as  $s_{i}=\eta_{i1}{\bf e}_{i}\overline{A_{0}}\, {}^t{\bf u}$.
We substitute~$ {\bf s}$ for (\ref{pdt}). Then, the right side of (\ref{pdt}) gives
\begin{gather}\label{13}
\big({\bf s}A\, {}^t{\bf t}\big)^N
=\big({\bf u}\zeta D_{2}\, {}^t{\bf t}\big)^N
=\big({\bf u}D_{2}\, {}^t{\bf t}\big)^N
=\sum_{\boldsymbol{m} \in X(n,N)}\binom{N}{\boldsymbol{m}}\prod_{i=0}^{n}(\eta_{i2} u_{i}t_{i})^{m_{i}}.
\end{gather}
Under this substitution, we consider the left side of (\ref{pdt}) and have
\begin{gather*}
\sum_{\boldsymbol{x} \in X(n,N)}\binom{N}{\boldsymbol{x}}\prod_{i=0}^{n-1}\big(s_{i}{\bf e}_{i} A_{0} \, {}^t{\bf t}\big)^{x_{i}}
=
\sum_{\boldsymbol{x} \in X(n,N)}\binom{N}{\boldsymbol{x}}
\prod_{i=0}^{n-1}s_{i}^{x_{i}}\prod_{i=0}^{n-1}\big({\bf e}_{i} A_{0} \, {}^t{\bf t}\big)^{x_{i}}\\
\hphantom{\sum_{\boldsymbol{x} \in X(n,N)}\binom{N}{\boldsymbol{x}}\prod_{i=0}^{n-1}\big(s_{i}{\bf e}_{i} A_{0} \, {}^t{\bf t}\big)^{x_{i}}}{}
=
\sum_{\boldsymbol{x} \in X(n,N)}\binom{N}{\boldsymbol{x}}
\prod_{i=0}^{n-1}\big(\eta_{i1}{\bf e}_{i}\overline{A_{0}}\, {}^t{\bf u}\big)^{x_{i}}\prod_{i=0}^{n-1}\big({\bf e}_{i} A_{0} \, {}^t{\bf t}\big)^{x_{i}}\\
\hphantom{\sum_{\boldsymbol{x} \in X(n,N)}\binom{N}{\boldsymbol{x}}\prod_{i=0}^{n-1}\big(s_{i}{\bf e}_{i} A_{0} \, {}^t{\bf t}\big)^{x_{i}}}{}
=
\sum_{\boldsymbol{x} \in X(n,N)}\binom{N}{\boldsymbol{x}}
\prod_{i=0}^{n-1}{\eta_{i1}}^{x_{i}}
\prod_{i=0}^{n-1}\big({\bf e}_{i}\overline{A_{0}}\, {}^t{\bf u}\big)^{x_{i}}\prod_{i=0}^{n-1}\big({\bf e}_{i} A_{0} \, {}^t{\bf t}\big)^{x_{i}}.
\end{gather*}
We expand the last two products of the above-mentioned
formula in terms of ${\bf u}^{\boldsymbol{m}}{\bf t}^{\boldsymbol{m}'} $'s;
\begin{gather*}
\prod_{i=0}^{n-1}\big({\bf e}_{i}\overline{A_{0}}\, {}^t{\bf u}\big)^{x_{i}}\prod_{i=0}^{n-1}\big({\bf e}_{i} A_{0} \, {}^t{\bf t}\big)^{x_{i}}
=\sum_{\boldsymbol{m},\boldsymbol{m}' \in X(n,N)}\binom{N}{\boldsymbol{m}}\binom{N}{\boldsymbol{m}'}\phi_{A}(\boldsymbol{x};
\boldsymbol{m})
\overline{\phi_{A}(\boldsymbol{x};\boldsymbol{m}')}{\bf u}^{\boldsymbol{m}}{\bf t}^{\boldsymbol{m}'}.
\end{gather*}
Now, the left side of (\ref{pdt}) gives
\begin{gather}\label{14}
\sum_{\boldsymbol{m},\boldsymbol{m}' \in X(n,N)}\binom{N}{\boldsymbol{m}}\binom{N}{\boldsymbol{m}'}\left( \sum_{\boldsymbol{x} \in X(n,N)} \prod_{i=0}^{n-1}\eta_{i1}^{x_{i}}\binom{N}{\boldsymbol{x}}\phi_{A}(\boldsymbol{x};\boldsymbol{m})
\overline{\phi_{A}(\boldsymbol{x};\boldsymbol{m}')}
\right)
{\bf u}^{\boldsymbol{m}}{\bf t}^{\boldsymbol{m}'}.
\end{gather}
By comparing coef\/f\/icients of ${\bf u}^k{\bf t}^{k'}$ of (\ref{13}) with (\ref{14}),
we conclude that
\begin{gather}\label{25}
\sum_{\boldsymbol{x} \in X(n,N)}
\prod_{i=0}^{n-1}
\eta_{i1}^{x_{i}}
\binom{N}{\boldsymbol{x}}
\phi_{A}(\boldsymbol{x};\boldsymbol{m})
\overline{\phi_{A}(\boldsymbol{x};\boldsymbol{m}')}
=
\frac{\prod\limits_{i=0}^{n-1}\eta_{i2}^{m_{i}}}
{\binom{N}{\boldsymbol{m}}}
\delta_{\boldsymbol{m} \boldsymbol{m}'}.
\end{gather}

Conversely, we assume  (\ref{25}) and substitute it for (\ref{14}).
Then, we reverse the above-mentioned computations and observe that
(\ref{pdt}) gives
\begin{gather*}
\big({\bf u}D_2 {\, {}^t\bf t}\big)^N
=\big({\bf u}A_{0}^{*}D_{1}A_{0}{\, {}^t\bf t}\big)^N.
\end{gather*}
This means
that
$
{\bf u}\zeta({\bf u},{\bf t}) D_2 {\, {}^t\bf t}
={\bf u}A_{0}^{*}D_{1}A_{0}{\, {}^t\bf t}
$
 holds for some $N$th root of unity $\zeta({\bf u},{\bf t})$.
 We have $A_{0}^{*}D_{1}A_{0}={\rm diag}(\zeta({\bf e}_{0},{\bf e}_{0}){\eta}_{02},\dots,\zeta({\bf e}_{n-1},{\bf e}_{n-1})
 {\eta}_{n-12})$.
We put ${\bf e}_{ij}(\theta,\varepsilon)=\cos{\theta}{\bf e}_{i}+ \varepsilon \sin{\theta}
{\bf e}_{j}$ for $|\varepsilon|>0$. We  def\/ine
\[
\zeta_{\varepsilon}(\theta)=
\frac
{{\bf e}_{ij}(\theta,\varepsilon)A_{0}^{*}D_{1}A_{0}  {\, {}^t{\bf e}_{ij}(\theta,\varepsilon)}}
{{\bf e}_{ij}(\theta,\varepsilon)D_2{\, {}^t{\bf e}_{ij}(\theta,\varepsilon)}
}=
\frac
{{\zeta({\bf e}_{i},{\bf e}_{i})\cos^2{\theta}}{\eta}_{i2}+\zeta({\bf e}_{j},{\bf e}_{j})\varepsilon^2\sin^2{\theta}{\eta}_{j2}}
{{\cos^2{\theta}}{\eta}_{i2}+\varepsilon^2\sin^2{\theta}{\eta}_{j2}}.
\]
Since $\zeta_{\varepsilon}(\frac{\pi}{2})=\zeta(\varepsilon{\bf e}_{j},\varepsilon{\bf e}_{j})$ does not depend on $\varepsilon$,
we have $\zeta(\varepsilon{\bf e}_{j},\varepsilon{\bf e}_{j})=\zeta({\bf e}_{j},{\bf e}_{j})$.
By taking $\varepsilon$ as ${\rm Arg}{(\frac{-\varepsilon^2 {\eta}_{j2}}{{\eta}_{i2}})}\not=0$, we have that
 $\zeta_{\varepsilon}(\theta)$ is a continuous function from $[0,\frac{\pi}{2}]$ to
the set of the $N$th root of unities.
Therefore $\zeta_{\varepsilon}(\theta)$ is a constant function, especially  $\zeta_{\varepsilon}(0)=
\zeta({\bf e}_{i},{\bf e}_{i})=\zeta_{\varepsilon}(\frac{\pi}{2})=\zeta(\varepsilon{\bf e}_{j},\varepsilon{\bf e}_{j})$
($0\leq i<j \leq n-1$).
Consequently we have
\begin{gather*}\zeta D_2
=A_{0}^{*}D_{1}A_{0}
\end{gather*}
for some $N$th root of unity $\zeta$.

\section{Proof of (\ref{hyp})}

Here we give a proof of the formula (\ref{hyp}) through direct computations. We need the following lemma.
\begin{lemma}\label{l3}
Put $\boldsymbol{p}=(p_{0},\dots,p_{n-1}) \in {\mathbb N}_{0}^{n}$
with
  $|\boldsymbol{p}|\leq N$.
For $\boldsymbol{m}=(m_{0},\dots,m_{n-1}) \in X(n,N)$
 and $\boldsymbol{z}=(z_{0},\dots,z_{n-1}) \in X(n,N-|\boldsymbol{p}|)$
 with $\boldsymbol{m}-\boldsymbol{z} \in {\mathbb N}_{0}^{n}$, we have
\begin{gather*}
\binom{N-|\boldsymbol{p}|}{\boldsymbol{z}}=
\binom{N}{\boldsymbol{m}}
\frac{\prod\limits_{i=0}^{n-1}(-m_{i})_{m_{i}-z_{i}}}{(-N)_{|\boldsymbol{p}|}}.
\end{gather*}
\end{lemma}

\begin{proof}
We compute
\begin{gather*}
\binom{N\!-\!|\boldsymbol{p}|}{\boldsymbol{z}} =
\binom{N}{\boldsymbol{m}}
\frac{(N-|\boldsymbol{p}|)!}{N!}
\prod_{i=0}^{n-1}\!\binom{m_{i}}{m_{i}\!-\!z_{i},z_{i}}
(m_{i}-z_{i})!
=
\binom{N}{\boldsymbol{m}}
\frac{\prod\limits_{i=0}^{n-1}\!(-1)^{m_{i}-z_{i}}(-m_{i})_{m_{i}-z_{i}}}{\binom{N}{|\boldsymbol{p}|}|\boldsymbol{p}|!}\\
\hphantom{\binom{N-|\boldsymbol{p}|}{\boldsymbol{z}}}{}
=
\binom{N}{\boldsymbol{m}}
\frac{(-1)^{|\boldsymbol{p}|}\prod\limits_{i=0}^{n-1}(-1)^{m_{i}-z_{i}}(-m_{i})_{m_{i}-z_{i}}}{(-N)_{|\boldsymbol{p}|}}
=
\binom{N}{\boldsymbol{m}}
\frac{\prod\limits_{i=0}^{n-1}(-m_{i})_{m_{i}-z_{i}}}{(-N)_{|\boldsymbol{p}|}}.\tag*{\qed}
\end{gather*}
\renewcommand{\qed}{}
\end{proof}
Now, we can prove the formula (\ref{hyp}).
We put $b_{ij}=1-a_{ij}$ and $\boldsymbol{c}_{i}=(c_{i0},\dots,c_{in-1})$. We compute
\begin{gather*}
\Phi_{N}(A;\boldsymbol{x}) =
\prod_{i=0}^{n-1}\left(\sum_{j=0}^{n-1}a_{ij}t_{j}\right)^{x_{i}}
=\prod_{i=0}^{n-1}
\left\{
\sum_{j=0}^{n-1}t_{j}
-
\sum_{j=0}^{n-1}b_{ij}t_{j}
\right\}^{x_{i}}\\
\phantom{\Phi_{N}(A;\boldsymbol{x})}
=\prod_{i=0}^{n-1}
\left\{
\sum_{p_{i}=0}^{x_{i}}
(-1)^{p_{i}}\binom{x_{i}}{x_{i}-p_{i},p_{i}}
\left(\sum_{j=0}^{n-1}t_{j}\right)^{x_{i}-p_{i}}
\left(\sum_{j=0}^{n-1}b_{ij}t_{j}\right)^{p_{i}}
\right\}\\
\phantom{\Phi_{N}(A;\boldsymbol{x})} =
\sum_{\substack{0 \leq p_{i}\leq x_{i}\\ (0 \leq i \leq n-1)}}
(-1)^{|\boldsymbol{p}|}
\left(\sum_{j=0}^{n-1}t_{j}\right)^{N-|p|}
\prod_{i=0}^{n-1}
\binom{x_{i}}{x_{i}-p_{i},p_{i}}
\left(\sum_{j=0}^{n-1}b_{ij}t_{j}\right)^{p_{i}}\\
\phantom{\Phi_{N}(A;\boldsymbol{x})} =
\sum_{\substack{0 \leq p_{i}\leq x_{i}\\ (0 \leq i \leq n-1)}}
(-1)^{|\boldsymbol{p}|}
\sum_{\substack{|\boldsymbol{c}_{i}|=p_{i}
,\\ |\boldsymbol{z}|=N-|p|}}
\binom{N-|\boldsymbol{p}|}{\boldsymbol{z}}{\boldsymbol{t}}^{\boldsymbol{z}}
\prod_{i=0}^{n-1}\binom{x_{i}}{x_{i}-p_{i},p_{i}}
\prod_{i=0}^{n-1}\binom{p_{i}}{\boldsymbol{c}_{i}}\prod_{i,j=0}^{n-1}b_{ij}^{c_{ij}}t_{j}^{c_{ij}}\\
\phantom{\Phi_{N}(A;\boldsymbol{x})}=
\sum_{0 \leq |\boldsymbol{p}|\leq N}
(-1)^{|\boldsymbol{p}|}
\sum_{\substack{|\boldsymbol{c}_{i}|=p_{i}
,\\ |\boldsymbol{z}|=N-|p|}}
\binom{N-|\boldsymbol{p}|}{\boldsymbol{z}}
\frac{\prod\limits_{i=0}^{n-1}(-1)^{|\boldsymbol{c}_{i}|}(-x_{i})_{|\boldsymbol{c}_{i}|}}{\prod\limits_{i,j=0}^{n-1}c_{ij}!}
{\boldsymbol{t}}^{\boldsymbol{z}}
\prod_{i,j=0}^{n-1}b_{ij}^{c_{ij}}t_{j}^{c_{ij}}\\
\phantom{\Phi_{N}(A;\boldsymbol{x})}
=\sum_{\boldsymbol{m} \in X(n,N)}
\binom{N}{\boldsymbol{m}}
\left[
\sum_{\sum_{{ij}}c_{ij}\leq N}
\frac{\prod\limits_{j=0}^{n-1}(-m_{j})_{\sum\limits_{i=0}^{n-1}c_{ij}}}{(-N)_{\sum_{{ij}}c_{ij}}}
\frac{\prod\limits_{i=1}^{n-1}(-x_{i})_{|\boldsymbol{c}_{i}|}}{\prod\limits_{i,j=0}^{n-1}c_{ij}!}\prod\limits_{i,j=0}^{n-1}b_{ij}^{c_{ij}}
\right]\boldsymbol{t}^{\boldsymbol{m}}.
\end{gather*}
In the last equation, we use
$|\boldsymbol{p}|=\sum_{ij}c_{ij}$, $m_{i}-z_{i}=\sum\limits_{i=0}^{n-1}c_{ij} \geq 0$
and Lemma~\ref{l3}.
Since $b_{0j}=b_{i0}=0$ for any $i$ and $j$, we have the formula.

\subsection*{Acknowledgements}
The author expresses his thanks to
two referees for their careful reading.
Without their many constructive comments and suggestions, this paper would not have been completed.

This work was supported by KAKENHI 21740032.

\pdfbookmark[1]{References}{ref}
\LastPageEnding


\begin{thebibliography}{99}

\footnotesize\itemsep=0pt

\bibitem{gr}
 Grif\/f\/iths R.C.,
Orthogonal polynomials on the multinomial distribution,
\href{http://dx.doi.org/10.1111/j.1467-842X.1971.tb01239.x}{{\it Austral. J. Statist.}} {\bf 13} (1971), 27--35.

\bibitem{rg}
 Gr\"{u}nbaum F.A., Rahman M.,
On a family of 2-variable orthogonal Krawtchouk polynomials,
\href{http://dx.doi.org/10.3842/SIGMA.2010.090}{{\it SIGMA}} {\bf 6} (2010), 090, 12~pages,
\href{http://arxiv.org/abs/1007.4327}{arXiv:1007.4327}.


\bibitem{m}
 Mizukawa H.,
 Zonal spherical functions on the  complex ref\/lection groups and $(m+1,n+1)$-hypergeometric functions,
\href{http://dx.doi.org/10.1016/S0001-8708(03)00092-6}{{\it Adv. Math.}} {\bf 184} (2004), 1--17.

 \bibitem{mt}
Mizukawa H., Tanaka H.,
 $(n+1,m+1)$-hypergeometric functions associated to character algebras,
\href{http://dx.doi.org/10.1090/S0002-9939-04-07399-X}{{\it Proc. Amer. Math. Soc.}} {\bf 132} (2004), 2613--2618.

\end{thebibliography}
\end{document}